\documentclass[12pt]{amsart}
\usepackage{amssymb,amsmath}
\usepackage{cite}

\newcommand{\dvg}{{\rm div}\,}
\newcommand{\N}{{\mathbb N}}

\newcommand{\R}{{\mathbb R}}

\newcommand{\eps}{\varepsilon}

\newcommand{\ws}[1]{\left|d#1\right|}

\makeatother

\title[Existence and symmetry of least energy solutions]{Existence and symmetry
of least energy solutions \\ for a class of quasi-linear elliptic equations}

\author{Louis Jeanjean}
\address{Laboratoire de Math\'ematiques (UMR 6623)
\newline\indent
Universit\'{e} de Franche-Comt\'{e}
\newline\indent
16, Route de Gray 25030 Besan\c{c}on Cedex, France}
\email{louis.jeanjean@univ-fcomte.fr}

\author{Marco Squassina}
\address{Dipartimento di Informatica
\newline\indent
Universit\`a degli Studi di Verona
\newline\indent
C\'a Vignal 2, Strada Le Grazie 15, 37134 Verona, Italy}
\email{marco.squassina@univr.it}

\thanks{The second author was partially supported by the
Italian PRIN Research Project 2005 {\em Metodi variazionali e topologici
nello studio di fenomeni non lineari}}

\newtheorem{theorem}{Theorem}[section]
\newtheorem{proposition}[theorem]{Proposition}

\newtheorem{lemma}[theorem]{Lemma}
\theoremstyle{definition}

\newtheorem{definition}[theorem]{Definition}

\newtheorem{remark}[theorem]{Remark}

\numberwithin{equation}{section}

\newcommand{\dom}[1]{{\rm dom}(#1)}
\newcommand{\gra}[1]{{\mathcal G}_{#1}}

\newcommand{\epi}[1]{{\rm epi}\left(#1\right)}

\newcommand{\dys}{\displaystyle}

\begin{document}



\begin{abstract}
Pour une large classe d'\'equations quasilin\'eaires elliptiques autonomes sur
$\R^N$, on montre l'existence d'une solution de moindre \'energie. On montre
aussi que toutes les solutions de moindres \'energies ont un signe constant et
sont, \`a une translation pr\`es, radiales. \vskip3pt \noindent For a general
class of autonomous quasi-linear elliptic equations on $\R^n$ we prove the
existence of a least energy solution and show that all least energy solutions
do not change sign and are radially symmetric up to a translation in $\R^n$.
\end{abstract}

\maketitle

\section{Introduction}
In this paper we show the existence, radial symmetry and sign of
the least energy solutions for a class of quasi-linear elliptic
equations,
\begin{equation}
\label{autonom} -\dvg(j_\xi(u,Du))+j_s(u,Du)=f(u)\quad \text{in
${\mathcal D}'(\R^n)$},
\end{equation}
where $\{\xi \mapsto j(s, \xi)\}$ is $p$-homogeneous. We look for
solutions of (\ref{autonom}) in $D^{1,p}(\R^n)$ where $1<p \leq
n$. If we set $F(s)=\int_0^sf(t)$, equation \eqref{autonom} is formally
associated with the functional
\begin{equation*}
I(u)=\int_{\R^n}j(u,Du)-\int_{\R^n}F(u),
\end{equation*}
which is nonsmooth on $D^{1,p}(\R^n)$, under natural growth
assumptions on the integrand $j(s,\xi)$ (see conditions
\eqref{growth1}-\eqref{growth2} below), although it admits
directional derivatives along the smooth directions. By least
energy solution of \eqref{autonom} we mean a nontrivial
function $u \in D^{1,p}(\R^n)$ such that
$$
I(u) = \inf\big\{ I(v) :\, \text{$v \in D^{1,p}(\R^n)$, $v \neq 0$
is a solution of \eqref{autonom}}\big\}.
$$
Our work is motivated by \cite{ByJeMa} where {\em abstract} conditions are
given under which problems of type (\ref{autonom}) admit a least
energy solution and all least energy solutions do not change sign
and are radially symmetric, up to a translation in $\R^n$. In this
paper, under quite general assumptions on $j(s,\xi)$ and $f(s)$,
we prove that, in fact, these {\em abstract} conditions hold.
In the special case of the $p$-Laplace equation ($1<p \leq n$)
\begin{equation}
\label{autonom1spLp} -\Delta_p u=f(u)\quad\text{in ${\mathcal
D}'(\R^n)$},
\end{equation}
there are various achievements regarding the existence of
solutions. For $p=2$, namely for
\begin{equation*}
-\Delta u=f(u)\quad\text{in ${\mathcal
D}'(\R^n)$},
\end{equation*}
we refer to the classical paper by Berestycki and Lions \cite{berlions} for the
scalar case and to the paper by Brezis and Lieb \cite{brezislieb} for both
scalar and systems cases.  In \cite{berlions,brezislieb} the existence of a
least energy solution is obtained. When $p \neq 2$ we refer to the papers
\cite{gazfer,gazsertan} and the references therein for the existence of
solutions. The issue of least energy solutions is not considered in these
papers. We also mention \cite{dome} where under assumptions on $f$,
allowing to work with regular functionals in $W^{1,p}(\R^n)$, the existence of
a least energy solution is derived. For a more general $j(s,\xi)$ the only
previous result about existence of least energy solutions is \cite[Theorem
3.2]{giacsqua} which actually played the r\^ole of a technical lemma therein.
However it requires significant restrictions on $f$ that we completely removed
in this paper.
\smallskip

In \cite{berlions,brezislieb} the existence of a least energy
solution is obtained by solving a constrained minimization problem
under suitable assumptions on $F$ and $f$. In \cite{brezislieb}
the authors assume that $F$ is a $C^1$ function on
$\R\setminus\{0\}$, locally Lipschitz around the origin and having
suitable sub-criticality controls at the origin and at infinity.
In Theorem 1 of \cite{gazfer}, the authors extend the existence
results of \cite{berlions} to the $p$-Laplacian case and need more
regularity on the function $F$ (for instance $f$ is taken in ${\rm
Lip}_{{\rm loc}}$). In our general setting, we consider a set of
assumptions on $F$ which is close to that of \cite{brezislieb} and
some natural assumptions on $j$ which are often considered in the
current literature of this kind of problems.

\subsection{Main result in the case $1<p<n$}
Let $F:\R\to\R$ be a function of class $C^1$ such that $F(0)=0$.
Denoting by $p^*$ the critical Sobolev exponent we assume that:
\begin{equation}
\label{zerocondit} \limsup_{s\to 0}\frac{F(s)}{|s|^{p^*}}\leq 0;
\end{equation}
\begin{equation}
\label{existposit} \text{there exists $s_0\in\R$ such that
$F(s_0)>0$.}
\end{equation}
Moreover, if $f(s)=F'(s)$ for any $s\in\R$,
\begin{equation}
\label{inftycondit} \lim_{s\to\infty}\frac{f(s)}{|s|^{p^*-1}}=0.
\end{equation}
Finally, \begin{equation} \label{manifold} \mbox{ if } u \in
D^{1,p}(\R^n) \mbox{ and } u \not \equiv 0 \mbox{ then } f(u) \not
\equiv 0.
\end{equation}
Condition (\ref{manifold}) is satisfied for instance if $f(s) \neq
0$ for $s \neq 0$ and small because if $u \in D^{1,p}(\R^n)$ and $
u \not \equiv 0$ the measure of the set $\{ x \in \R^n: \eta \leq
|f(u(x))| \leq 2 \eta\}$ is positive for $\eta >0$ and small.
\smallskip

Let $j(s,\xi):\R\times\R^n\to\R$ be a function of class $C^1$ in
$s$ and $\xi$ and denote by $j_s$ and $j_\xi$ the derivatives of
$j$ with respect of $s$ and $\xi$ respectively. We assume that:
\begin{equation}
\label{stricconv} \text{for all $s\in\R$ the map $\{\xi\mapsto
j(s,\xi)\}$ is strictly convex and $p$-homogeneous;}
\end{equation}
there exist positive constants $c_1,c_2,c_3,c_4$ and $R$ such that
\begin{equation}
\label{growth1} c_1|\xi|^p\leq j(s,\xi)\leq
c_2|\xi|^p,\qquad\text{for all $s\in\R$ and $\xi\in\R^n$};
\end{equation}
\begin{equation}
\label{growth2} |j_s(s,\xi)|\leq c_3|\xi|^p,\quad
|j_\xi(s,\xi)|\leq c_4|\xi|^{p-1}, \qquad\text{for all $s\in\R$
and $\xi\in\R^n$};
\end{equation}
\begin{equation}
\label{posit} j_s(s,\xi)s\geq 0,\qquad\text{for all $s\in\R$ with
$|s|\geq R$ and $\xi\in\R^n$}.
\end{equation}
Conditions \eqref{stricconv}-\eqref{posit} on $j$ are quite
natural assumptions and were already used, e.g., in
\cite{ab1,canino,cd,squastop,squastoul}. See also
Remark~\ref{rem27} for further comments regarding
the r\^ole played by condition~\eqref{posit}. \vskip5pt

\begin{theorem}
\label{MainThm} Assume that conditions
\eqref{zerocondit}-\eqref{posit} hold. Then equation
\eqref{autonom} admits a least energy solution $u\in
D^{1,p}(\R^n)$. Furthermore any least energy solution of
\eqref{autonom} has a constant sign and is radially symmetric, up
to a translation in $\R^n$.
\end{theorem}

\subsection{Main result in the case $p=n$}

Let $F:\R\to\R$ be a function of class $C^1$ such that $F(0)=0$.
We assume that:
\begin{equation}
\label{zerocondit2} \text{there exists $\delta>0$ such that
$F(s)<0$ for all $0<|s|\leq\delta$};
\end{equation}
\begin{equation}
\label{existposit2} \text{there exists $s_0\in\R$ such that
$F(s_0)>0$};
\end{equation}
\begin{equation}
\label{inftycondit2} \text{there exist $q>1$ and $c>0$ such that
$|f(s)|\leq c+c|s|^{q-1}$ for all $s\in\R$}.
\end{equation}
\begin{equation} \label{manifold2} \mbox{ if } u \in
D^{1,n}(\R^n) \mbox{ and } u \not \equiv 0 \mbox{ then } f(u) \not
\equiv 0.
\end{equation}
Concerning the Lagrangian $j$ we still assume conditions
\eqref{stricconv}-\eqref{posit} (with $p=n$). \vskip3pt

\begin{theorem}
\label{MainThm2} Assume that conditions
\eqref{stricconv}-\eqref{manifold2} hold. Then equation
\eqref{autonom} admits a least energy solution $u\in
D^{1,n}(\R^n)$. Furthermore any least energy solution of
\eqref{autonom} has a constant sign and if a least energy solution
$u \in D^{1,n}(\R^n)$ satisfies $u(x) \to 0$ as $|x| \to \infty$ then
it is radially symmetric, up to a translation in $\R^n$.
\end{theorem}

Our approach to prove Theorems \ref{MainThm} and \ref{MainThm2} is
based in an essential way on the work \cite{ByJeMa}. There {\it
abstract} conditions (see {\bf (C1)-(C3)} and {\bf (D1)-(D3)}
below) are given which, if they are satisfied, guarantee the
conclusions of our Theorems \ref{MainThm} and \ref{MainThm2}.  \medskip

We point out that the way we prove the existence of least energy
solutions, by solving a constrained minimization problem, is
crucial in order to get the symmetry and sign results of all least
energy solutions. First we show that the problem
\begin{equation}
\label{minconstprob} \min\left\{\int_{\R^n}j(u,Du):\,u\in
D^{1,p}(\R^n),\,\, F(u)\in L^1(\R^n),\,\,\int_{\R^n}F(u)=
1\right\}
\end{equation}
admits a solution, which is the hardest step. To do this, we exploit some tools from
non-smooth critical point theory, such as the weak slope, developed
in \cite{cdm,dm,ioffe,katriel} (see Section \ref{nonsmoothsect}).   Then
we prove that any minimizer is of class $C^1$ and satisfy the
Euler-Lagrange equation as well as the Pucci-Serrin identity. This
allow us to check the {\em abstract} conditions of \cite{ByJeMa}
which provide a link between least action solutions of
\eqref{autonom} and  solutions of problem \eqref{minconstprob}.
Roughly speaking if the {\em abstract} conditions hold then there
exist a least energy solution and to any least energy solution of
\eqref{autonom} correspond, up to a rescaling, a minimizer of
\eqref{minconstprob}. It is proved in \cite{Ma} that any such
minimizer are radially symmetric. In addition it is shown in
\cite{ByJeMa} that any minimizer has a constant sign.
\medskip

Let us point out that, in our setting, the existence results that
we obtain have no equivalent in the literature. Also, even
assuming the existence of least energy solutions, our results of
symmetry and sign are new. In particular we observe that, under
our assumptions, to try to show that they are radial using moving
plane methods or rearrangements arguments is hopeless. We
definitely need to use the approach of \cite{ByJeMa} which, in
turn, is based on the remarkable paper \cite{Ma}. In \cite{Ma}
results of symmetry for $C^1$ mimimizers are obtained for general
functionals under one or several constraints. Let us finally
mention that in \cite{ByJeMa}, and thus in our paper, the results
of radial symmetry are obtained without using the fact that our
solutions have a constant sign.

\section{The case $1<p<n$}

\subsection{Conditions $\boldsymbol{(C1)}$-$\boldsymbol{(C3)}$}
\noindent
In the following $D^{1,p}(\R^n)$ will denote the closure of the space
$C^\infty_c(\R^n)$ with respect to the norm $\|u\|=(\int_{\R^n}|Du|^p)^{1/p}$
and $D^*$ is the dual space of $D^{1,p}(\R^n)$. Let us consider the problem
\begin{equation}
\label{autonomG}
-\dvg(j_\xi(u,Du))+j_s(u,Du)=f(u)\quad
\text{in ${\mathcal D}'(\R^n)$},
\end{equation}
associated with the functional
$$
I(u)=\int_{\R^n}j(u,Du)-\int_{\R^n}F(u),
$$
where $F(s)=\int_0^s f(t)dt$. Moreover, introducing the functionals,
\begin{equation*}
J(u)=\int_{\R^n}j(u,Du),\quad
V(u)=\int_{\R^n}F(u),\qquad u\in D^{1,p}(\R^n),
\end{equation*}
we  consider the following constrained problem
\begin{equation}
\label{constminim}
\tag{$P_1$}
\text{minimize $J(u)$ subject to the constraint $V(u)=1$}.
\end{equation}
More precisely, let us set
\begin{equation*}
X=\big\{u\in D^{1,p}(\R^n):\, F(u)\in L^1(\R^n)\big\},
\end{equation*}
and
$$
T=\inf_{{\mathcal C}}J,\qquad {\mathcal C}=\big\{u\in X:\,V(u)=1\big\}.
$$
Consider the following conditions:
\vskip6pt
\noindent
$\boldsymbol{(C1)}$ $T>0$ and problem $(P_1)$ has a minimizer $u\in X$;
\vskip8pt
\noindent
$\boldsymbol{(C2)}$ any minimizer $u\in X$ of $(P_1)$ is a $C^1$ solution and
satisfies the equation
\begin{equation}
\label{autonomGC}
-\dvg(j_\xi(u,Du))+j_s(u,Du)=\mu f(u)\quad
\text{in ${\mathcal D}'(\R^n)$},
\end{equation}
for some $\mu\in\R$.
\vskip8pt
\noindent
$\boldsymbol{(C3)}$ any solution $u\in X$ of equation \eqref{autonomGC} satisfies the identity
$$
(n-p)J(u)=\mu nV(u).
$$
\vskip2pt \noindent

From \cite{ByJeMa} we have the following

\begin{proposition} \label{abstract}
Assume that $1<p<n$ and that conditions {\bf (C1)-(C3)} hold. Then
(\ref{autonom}) admits a least energy solution and each least
energy solution has a constant sign and is radially symmetric, up
to a translation in $\R^n$.
\end{proposition}

Indeed $X$ is an admissible function space in the
sense introduced in \cite{ByJeMa}. Then Proposition 3 of \cite{ByJeMa} gives
the existence of a least energy solution and that any least energy
solution is radially symmetric. Finally, the sign result follows
directly from Proposition 5 of \cite{ByJeMa}. \medskip

In view of Proposition \ref{abstract}, our aim is now to prove that
conditions {\bf (C1)-(C3)} are fulfilled under assumptions
\eqref{zerocondit}-\eqref{posit}.

\begin{remark} \label{adaptation}
In \cite{ByJeMa}, in the scalar case, the equation considered is precisely
\eqref{autonom1spLp} and the corresponding functional
$$
I(u) = \frac{1}{p}\int_{\R^n}|\nabla u|^p - \int_{\R^n}F(u).
$$
However, in order to show that conditions {\bf (C1)-(C3)} imply the
conclusion of Proposition \ref{abstract}, the only property of
$|\nabla u|^p$ that it is used is $p$-homogeneity, namely
that \eqref{stricconv} hold.
\end{remark}

\smallskip

\subsection{Some recalls of non-smooth critical point theory}
\label{nonsmoothsect}
In this section we recall some abstract notions that will be used in the sequel.
We refer the reader to \cite{cdm,dm,ioffe,katriel}, where this theory is fully developed.
\vskip2pt
Let $X$ be a metric space and let $f:X\to\R\cup\{+\infty\}$
be a lower semicontinuous function. We set
\begin{equation*}
\dom{f} = \left\{u\in X:\,f(u)<+\infty\right\}
\quad\text{and}\quad
\epi{f}=\left\{(u,\eta)\in X\times\R:\,f(u)\leq\eta\right\}.
\end{equation*}
The set $\epi{f}$ is endowed with the metric
$$
d\left((u,\eta),(v,\mu)\right)=
\left(d(u,v)^2+(\eta-\mu)^2\right)^{1/2}.
$$
Let us define the function $\gra{f}:\epi{f}\to\R$ by setting
\begin{equation*}
\gra{f}(u,\eta)=\eta.
\end{equation*}
Note that $\gra{f}$ is Lipschitz  continuous of constant $1$.
In the following $B(u,\delta)$ denotes the open ball of center
$u$ and of radius $\delta$. We recall the definition of the weak slope for a
continuous function.

\begin{definition}\label{defslope}
Let $X$ be a complete metric space, $g:X \to \R $ a continuous function,
and $u\in X$. We denote by
$|dg|(u)$ the supremum of the real numbers $ \sigma$ in
$[0,\infty)$ such that there exist $ \delta >0$
and a continuous map
$$
{\mathcal H}\,:\,B(u, \delta) \times[ 0, \delta]  \to X,
$$
such that, for every $v$ in $B(u,\delta) $, and for every
$t$ in $[0,\delta]$ it results
\begin{align*}
d({\mathcal H}(v,t),v) &\leq t,\\
g({\mathcal H}(v,t)) &\leq g(v)-\sigma t.
\end{align*}
The extended real number $|dg|(u)$ is called the weak
slope of $g$ at $u$.
\end{definition}

According to the previous definition, for every lower
semicontinuous function $f$ we can consider the metric
space $\epi{f}$ so that the weak slope of $\gra{f}$ is
well defined. Therefore, we can define the weak slope of a  lower semicontinuous function $f$
by using $|d\gra{f}|(u,f(u))$.

\begin{definition}
\label{defnwslsc}
For every $u\in\dom{f}$ let
\begin{equation*}
\ws{f}(u)=
\begin{cases}
\dys \frac{\ws{\gra{f}}(u,f(u))}{{\sqrt{1-\ws{\gra{f}}(u,f(u))^2}}},
& \text{if $\ws{\gra{f}}(u,f(u))<1$}, \\
\noalign{\vskip4pt}
+\infty, & \text{if $\ws{\gra{f}}(u,f(u))=1$}.
\end{cases}
\end{equation*}
\end{definition}

The previous notion allows to give, in this framework, the definition of critical point of $f$
(namely a point $u\in\dom{f}$ with $\ws{f}(u)=0$) as well as the following

\begin{definition}\label{defps}
Let $X$ be a complete metric space, $f:X\to \R\cup\{+\infty\}$ a lower semicontinuous function and
let $c\in\R$. We say that $f$ satisfies the Palais-Smale condition at level $c$ ($(PS)_c$
in short), if every sequence $(u_n)$ in $\dom{f}$ such that $\ws{f}(u_n)\to 0$ and $f(u_n)\to c$
($(PS)_c$ sequence, in short) admits a subsequence $(u_{n_k})$ converging in $X$.
\end{definition}

We now recall a consequence of Ekeland's variational principle \cite{EVP} in the framework
of the weak slope (just apply \cite[Theorem 3.3]{dm} with $r=r_h=\sigma=\sigma_h=\eps_h$ for a sequence
$\eps_h\to 0$; see also \cite[Corollary 3.4]{dm}).

\begin{proposition}
\label{ekelVP} Let $X$ be a complete metric space and $f:X\to
\R\cup\{+\infty\}$ a lower semicontinuous function which is
bounded from below. Assume that $(u_h)\subset\dom{f}$ is a
minimizing sequence for $f$, that is $f(u_h)\to c=\inf_Xf$. Then
there exists a sequence $(\eps_h)\subset\R^+$ with $\eps_h\to 0$
as $h\to\infty$ and a sequence $(v_h)\subset X$ such that
$$
\ws{f}(v_h)\leq\eps_h,\qquad d(v_h,u_h)\leq\eps_h,\qquad f(v_h)\leq f(u_h).
$$
In particular $(v_h)$ is a minimizing sequence and a $(PS)_c$ sequence for $f$.
\end{proposition}

Finally, we mention the notion of subdifferential as introduced in \cite{campa}.

\begin{definition}
For a function $f:X\to\R$, we set
$$
\partial f(x)=\big\{\alpha\in X': (\alpha,-1)\in N_{{\rm epi}(f)}(x,f(x))\big\},
$$
where $X'$ is the dual space of $X$,
$N_C(x)=\{\nu\in X':\langle\nu,v\rangle\leq 0,\,\,\text{for all $v\in T_C(x)$}\}$
is the normal cone (and $T_C(x)$ the tangent cone) to the set $C$ at the point $x$.
\end{definition}

More precisely, see \cite[Definition 3.1, 3.3 and 4.1]{campa}.
\medskip

\subsection{Verification of conditions $\boldsymbol{(C1)}$-$\boldsymbol{(C3)}$}
In this section we assume that \eqref{zerocondit}-\eqref{posit} hold and then show
that $\boldsymbol{(C1)}$-$\boldsymbol{(C3)}$ are fulfilled.
\medskip

First we extend $J{|_{\mathcal C}}$ to the functional $J^*:D^{1,p}(\R^n)\to\R\cup\{+\infty\}$,
\begin{equation}
\label{extension}
J^*(u)=
\begin{cases}
J(u) & \text{if\enskip $u\in {\mathcal C}$} \\
+\infty  & \text{if\enskip $u\not\in {\mathcal C},$}
\end{cases}
\end{equation}
which turns out to be lower semicontinuous.
Since ${\mathcal C}$ can be regarded as a metric space endowed with the metric of $D^{1,p}(\R^n)$,
the weak slope $\ws{J{|_{\mathcal C}}}(u)$ and the Palais-Smale condition for $J{|_{\mathcal C}}$ may be defined.

\begin{lemma}
\label{quasieuler}
For all $u\in {\mathcal C}$ there exists $\mu\in\R$ such that
$$
\ws{J{|_{\mathcal C}}}(u)\geq\sup\left\{J'(u)(v)-\mu V'(u)(v):\,v\in C^{\infty}_c(\R^n),\,\|Dv\|_{p}\leq 1\right\}\,.
$$
In particular, for each $(PS)_c$-sequence $(u_h)$ for $J{|_{\mathcal C}}$ there exists $(\mu_h)\subset\R$ such that
$$
\lim_h\sup\left\{J'(u_h)(v)-\mu_h V'(u_h)(v):\,v\in C^{\infty}_c(\R^n),\,\|Dv\|_{p}\leq 1\right\}=0.
$$
\end{lemma}
\begin{proof}
By condition~\eqref{growth2}, for all $u\in {\mathcal C}$ and any
$v\in C^{\infty}_c(\R^n)$ the directional derivative of $J$ at $u$
along $v$ exists and it is given by
\begin{equation}
\label{gradient} J'(u)(v)=\int_{\R^n}j_{\xi}(u,Du)\cdot
Dv+\int_{\R^n}j_s(u, Du)v.
\end{equation}
Moreover the function $\{u\mapsto J'(u)(v)\}$ is continuous from ${\mathcal C}$ into $\R$.
Of course, we may assume that $|dJ{|_{\mathcal C}}|(u)<+\infty.$ If $J^*$ is defined as in \eqref{extension}, we have $|dJ^*|(u)=|dJ{|_{\mathcal C}}|(u)$,
so that by virtue of  \cite[Theorem 4.13]{campa} there exists $\omega\in\partial J^*(u)$ with $|dJ^*|(u)\geq\|\omega\|_{D^*}.$
Moreover, by \cite[Corollary 5.9(ii)]{campa}, we have $\partial J^*(u)\subseteq \partial J(u)+\R V'(u)$.
Finally, by \cite[Theorem 6.1(ii)]{campa}, we get $\partial J(u)=\{\eta\}$ where, for any function $v\in C^\infty_c(\R^n)$,
$\langle\eta,v\rangle=J'(u)(v)$. This concludes the proof.
\end{proof}

\begin{remark}
\label{modifEq} Assuming only that $F$ is $C^1$ on
$\R\setminus\{0\}$ and it is locally Lipschitz around the origin
(as in \cite{brezislieb} as it follows by \cite[assumption
2.8]{brezislieb} which is used in the proof of Theorem 2.2
therein)  Lemma \ref{quasieuler} cannot hold in the form it is
stated. In this more general case, there would exist $\mu\in\R$
and some function $\varphi\in L^\infty(\R^n)$ such that the
solutions of the minimum problem satisfy
\begin{equation}
\label{mgc}
-\dvg(j_\xi(u,Du))+j_s(u,Du)=\mu f(u)\chi_{\{u\neq 0\}}+\varphi\chi_{\{u=0\}}\quad
\text{in ${\mathcal D}'(\R^n)$}.
\end{equation}
In fact, notice that this is exactly what is obtained at the
bottom of page 103 in \cite{brezislieb}. Then, in light of the strong
regularity of their solutions (that is $W^{2,\sigma}_{{\rm loc}}(\R^n)$ for some
$\sigma>1$), the equation is satisfied pointwise and as $\Delta
u=0$ a.e. in $\{u=0\}$ (by a result of Stampacchia, see \cite{stampac})
they infer $\varphi=0$. On the other hand, in our degenerate
framework we cannot reach this regularity level and concluding
that $\varphi=0$ (and hence that $u$ solves \eqref{autonomGC})
seems, so far, out of reach.
\end{remark}

Now we recall (see \cite[Theorem 2]{squastop}) the following
\begin{lemma}
\label{compscalar}
Let $(u_h)$ be a bounded sequence in $D^{1,p}(\R^n)$ and, for each $v\in C^\infty_c(\R^n)$, set
\begin{equation*}
\langle w_h,v\rangle=\int_{\R^n}j_\xi(u_h,Du_h)\cdot Dv+\int_{\R^n} j_s(u_h, Du_h)v=J'(u_h)(v).
\end{equation*}
If the sequence $(w_h)$ is strongly convergent to some $w$ in $D^*(\Omega)$ for
each open and bounded subset $\Omega\subset\R^n$, then $(u_h)$ admits a
strongly convergent subsequence in $D^{1,p}(\Omega)$ for each open
and bounded subset $\Omega\subset\R^n$.
\end{lemma}

\begin{lemma}
\label{lemma1}
Assume \eqref{zerocondit}-\eqref{posit}. Then condition $\boldsymbol{(C1)}$ holds.
\end{lemma}
\begin{proof}
In view of assumption \eqref{existposit} the constraint ${\mathcal C}$ is not empty (see Step 1 at page 324 in \cite{berlions}).
Let then $(u_h)\subset {\mathcal C}$ be a minimizing sequence for $J|_{{\mathcal C}}$. Therefore, we have
$$
\lim_{h}\int_{\R^n}j(u_h,Du_h)=T,\quad
F(u_h)\in L^1(\R^n),\quad
\int_{\R^n}F(u_h)=1,\quad\text{for all $h\in\N$}.
$$
After extracting a subsequence, still denoted by $(u_h)$, we get
using \eqref{growth1},
\begin{equation}
\label{subseq}
u_h\rightharpoonup u\quad \text{in $L^{p^*}(\R^n)$},\quad
Du_h\rightharpoonup Du\quad \text{in $L^{p}(\R^n)$},\quad
u_h(x)\to u(x)\quad\text{a.e.}
\end{equation}
As $j(s,\xi)$ is positive, convex in the $\xi$ argument,
$u_h\to u$ in $L^1_{{\rm loc}}(\R^n)$ and $Du_h\rightharpoonup Du$ in $L^1_{{\rm loc}}(\R^n)$,
by well known lower semicontinuity results (cf.\ \cite{ioffe1,ioffe2}), it follows
\begin{equation}
\label{lowersem}
\int_{\R^n}j(u,Du)\leq \liminf_{h}\int_{\R^n}j(u_h,Du_h)=T.
\end{equation}
Moreover, setting $F=F^+-F^-$ with $F^+ = \max\{F,0\}$ and $F^- =
\max\{-F,0\}$, in view of assumptions \eqref{zerocondit} and
\eqref{inftycondit} (which implies that $|F(s)|/|s|^{p^*}$ goes to
zero as $s\to\infty$), fixing some $c>0$ one can find $r_2>r_1>0$
such that $F^+(s)\leq c|s|^{p^*}$ for all $|s|\leq r_1$ and
$|s|\geq r_2$, so that
$$
1+\int_{\R^n}F^-(u_h)=\int_{\R^n}F^+(u_h)\leq c\int_{\R^n}|u_h|^{p^*}\chi_{\{|u_h|\leq r_1\}\cup\{|u_h|\geq r_2\}}+
\beta{\mathcal L}^n(\{|u_h|>r_1 \})
$$
where $\beta=\max\{F^+(s):r_1\leq |s|\leq r_2\}$ and ${\mathcal
L}^n$ is the Lebesgue measure in $\R^n$. Clearly ${\mathcal
L}^n(\{|u_h|>r_1 \})$ remains uniformly bounded, as $(u_h)$ is
bounded in $L^{p^*}(\R^n)$. Hence, by Fatou's lemma, this yields
$F^+(u),F^-(u)\in L^1(\R^n)$ and thus, finally, $F(u)\in
L^1(\R^n)$. We have proved that $u\in X$. Notice that, still by
assumptions \eqref{zerocondit} and \eqref{inftycondit}, in light
of \cite[Lemma 2.1]{brezislieb}, we find two positive constants
$\eps_1,\eps_2$ such that
\begin{equation}
\label{cond-misura}
{\mathcal L}^n\big(\{x\in\R^n:|u_h(x)|>\eps_1 \}  \big)\geq \eps_2,\qquad\text{for all $h\in\N$}.
\end{equation}
Hence, in view of \cite[Lemma 6]{lieb-inv} (cf.\ the proof due to
H.\ Brezis at the end of page 447 in \cite{lieb-inv}) there exists
a shifting sequence $(\xi_h)\subset\R^n$ such that
$(u_h(x+\xi_h))$ converges weakly to a nontrivial limit. Thus, in
\eqref{subseq}, we may assume that $u\not\equiv 0$. Applying
Proposition \ref{ekelVP} to the lower semicontinuous functional
$J^*$ defined in  \eqref{extension}, we can replace the minimizing
sequence $(u_h)\subset {\mathcal C}$ by a minimizing sequence
$(v_h)\subset {\mathcal C}$  with $\|v_h-u_h\|_{D^{1,p}}=o(1)$ as
$h\to\infty$ (we shall rename $v_h$ again as $u_h$) such that the
weak slope vanishes, namely $|dJ|_{{\mathcal C}}|(u_h)\leq\eps_h$,
with $\eps_h\to 0$ as $h\to\infty$. It follows by
Lemma~\ref{quasieuler} that there exists a sequence
$(\mu_h)\subset\R$ of Lagrange multipliers such that
\begin{equation}
\label{quasi-sol}
J'(u_h)(v)=\mu_h V'(u_h)(v)+\langle \eta_h,v\rangle,\qquad\text{for all $h\in\N$ and $v\in C^{\infty}_c(\R^n)$,}
\end{equation}
where $\eta_h$ strongly converges to $0$ in $D^*$. Also for any
bounded domain $\Omega \subset \R^n$, by \eqref{inftycondit}
(which implies that, for each $\eps>0$, there exists $a_\eps\in\R$
such that $|f(s)|\leq a_\eps+\eps|s|^{p^*-1}$ for all $s\in\R$) it
follows that the map $D^{1,p}(\Omega)\ni v\mapsto f(v)\in
D^*(\Omega)$ is completely continuous. Thus by condition \eqref{manifold},
since $u \neq 0$, there exists a function $\psi_0\in
C^\infty_c(\R^n)$ such that, setting $K_0={\rm supp}(\psi_0)$, it
holds
$$
V'(u_h)(\psi_0)=\int_{K_0} f(u_h)\psi_0\not\to 0,\quad\text{as
$h\to\infty$}.
$$
Also the sequence $(J'(u_h)(\psi_0))$ is bounded. In fact, denoting by $C$
a generic positive constant, we have by \eqref{growth2} that
\begin{gather*}
\int_{K_0} |j_\xi(u_h,Du_h)||D\psi_0|\leq C\int_{K_0} |Du_h|^{p-1}\leq C\Big(\int_{K_0} |Du_h|^{p}\Big)^{\frac{p-1}{p}}\leq C,  \\
\int_{K_0} |j_s(u_h,Du_h)| |\psi_0| \leq C \int_{K_0} |Du_h|^p \leq C,
\end{gather*}
and we conclude using \eqref{gradient}. Now, formula \eqref{quasi-sol} yields
$$
\mu_h V'(u_h)(\psi_0)+\langle \eta_h,\psi_0\rangle=J'(u_h)(\psi_0)
$$
with $\eta_h\to 0$ in $D^*$. We deduce that the sequence $(\mu_h)$
is bounded in $\R$ and thus we can assume that it converges to
some $\mu\in\R$. It follows that $w_h=\mu_h V'(u_h)+\eta_h=\mu_h
f(u_h)+\eta_h$ converges strongly to some $w$ in $D^*(\Omega)$.
Therefore, by Lemma~\ref{compscalar}, we infer that $(u_h)$
admits a subsequence which strongly converges in
$D^{1,p}(\Omega)$. Thus, we have proved that the sequence $(u_h)$
is locally compact in $D^{1,p}(\R^n)$ (fact that is useful in the forthcoming steps).  From this we can easily
deduce that $J'(u)(v) = \mu V'(u)(v)$ for all $v \in
C_c^{\infty}(\R^n)$, namely \eqref{autonomGC}. However this is not
enough (nor necessary) to show that $u$ is  a minimizer of $(P_1)$
since we do not know if $V(u)=1$.  In this aim  we set
$u_\sigma(x)=u(x/\sigma)$. Then it holds $J(u_\sigma)=\sigma^{n-p}J(u)$
and $V(u_\sigma)=\sigma^{n}V(u)$ and hence, by a simple scaling
argument, we get that
\begin{equation}
\label{rescdis} \int_{\R^n}j(w,Dw)\geq
T\left(\int_{\R^n}F(w)\right)^{\frac{n-p}{n}}\!\!\!\!,\qquad
\text{for all $w\in D^{1,p}(\R^n)$ with $V(w)>0$}.
\end{equation}
We follow now an argument in the spirit of the perturbation method
developed in the proof of \cite[Lemma 2.3]{brezislieb}. Taking any
function $\phi\in L^{p^*}(\R^n)$ with compact support, we claim
that
\begin{equation*}
\int_{\R^n}F(u_h+\phi)=1+\int_{\R^n}F(u+\phi)-\int_{\R^n}F(u)+o(1),\quad\text{as $h\to\infty$}.
\end{equation*}
Indeed, if we set $K={\rm supp}(\phi)$, we get
\begin{align*}
\int_{\R^n}F(u_h+\phi) &=\int_{\R^n}F(u_h)+\int_{K}F(u_h+\phi)-F(u_h) \\
&=1+\int_{K}F(u+\phi)-F(u)+o(1) \\
&=1+\int_{\R^n}F(u+\phi)-F(u)+o(1),
\quad\text{as $h\to\infty$},
\end{align*}
where the second equality follows by the dominated convergence theorem
in light of \eqref{inftycondit} and the strong convergence of $u_h$ to $u$ in $L^{p^*-1}(K)$. Moreover, we have
\begin{align*}
\int_{\R^n}j(u_h+\phi,Du_h+D\phi)&=T+\int_{\R^n}j(u_h+\phi,Du_h+D\phi)-j(u_h,Du_h)+o(1) \\
&=T+\int_{K}j(u_h+\phi,Du_h+D\phi)-j(u_h,Du_h)+o(1) \\
&=T+\int_{K}j(u+\phi,Du+D\phi)-j(u,Du)+o(1) \\
&=T+\int_{\R^n}j(u+\phi,Du+D\phi)-j(u,Du)+o(1),
\end{align*}
as $h\to\infty$, where the third equality is justified again by
the dominated convergence theorem, since as $Du_h\to Du$ in
$L^p(K)$ for $h\to\infty$ we have
\begin{gather*}
|j(u_h+\phi,Du_h+D\phi)-j(u_h,Du_h)|\leq c_2|Du_h+D\phi|^p\chi_{K}+c_2|Du_h|^p\chi_{K}, \\
j(u_h+\phi,Du_h+D\phi)-j(u_h,Du_h)\to j(u+\phi,Du+D\phi)-j(u,Du)\quad\text{a.e.\ in $\R^n$}.
\end{gather*}
Therefore, choosing $w=u_h+\phi$ inside inequality \eqref{rescdis}, where $\phi\in D^{1,p}(\R^n)$
has compact support and $1+\int_{\R^n}F(u+\phi)-F(u)>0$, and taking the limit as $h\to\infty$, it follows that
\begin{equation*}
T+\int_{\R^n}j(u+\phi,Du+D\phi)-\int_{\R^n} j(u,Du)\geq T\left( 1+\int_{\R^n}F(u+\phi)-\int_{\R^n}F(u)  \right)^{\frac{n-p}{n}}.
\end{equation*}
Fixed $\lambda$ close to $1$, we consider for some $r>1$ a $C^\infty$ function $\Lambda:\R^+\to\R^+$ such that
\begin{equation}
\label{defLambda} \text{$\Lambda(t)=\lambda$\,\,\, if $t\leq
1$,\quad $\Lambda(t)=1$\,\,\, if $t\geq r$,\quad
$\rho=\inf_{t\in\R^+}\Lambda(t)>\frac{1}{2}$,\quad $\sup_{t \in\R^+}|\Lambda'(t)| < \frac{\rho}{r}$,}
\end{equation}
and we introduce the smooth and bijective map $\Pi:\R^n\to\R^n$ by setting $\Pi(x)=\Lambda(|x|)x$,
for all $x\in\R^n$. Finally, we set
$$
\Pi_h(x)=h\Pi\left(\frac{x}{h}\right)=
\begin{cases}
\lambda x & \text{if $|x|\leq h$,} \\
\Lambda\big(\frac{|x|}{h}\big)x  & \text{if $h\leq |x|\leq rh$,} \\
x & \text{if $|x|\geq rh$,} \\
\end{cases}
\qquad
\phi_h(x)=u(\Pi_h(x))-u(x).
$$
In particular, it follows that $\phi_h\in D^{1,p}(\R^n)$ is a
compact support function which satisfies $1+\int_{\R^n}F(u+\phi_h)-F(u)>0$, at least for all values of
$\lambda$ sufficiently close to $1$ (see equation \eqref{Fasssintt} below). Hence, for any $h\in\N$, we conclude
\begin{equation*}
T+\int_{\R^n}j(u+\phi_h,Du+D\phi_h)-\int_{\R^n} j(u,Du)\geq T\left( 1+\int_{\R^n}F(u+\phi_h)-\int_{\R^n}F(u)  \right)^{\frac{n-p}{n}}\!\!\!\!.
\end{equation*}
Notice that, we have
\begin{equation*}
\int_{\R^n}j(u+\phi_h,Du+D\phi_h)=\int_{\R^n}j(u(\Pi_h(x)),Du(\Pi_h(x)))=I_1+I_2+I_2.
\end{equation*}
In view of assumptions \eqref{stricconv} and \eqref{growth1}, by dominated convergence we have
\begin{align*}
I_1 &=\int_{\R^n}j(u(\lambda x),\lambda(Du)(\lambda x))\chi_{\{|x|\leq h\}} \\
&=\lambda^p\int_{\R^n}j(u(\lambda x),(Du)(\lambda x))\chi_{\{|x|\leq h\}} \\
&=\lambda^p\int_{\R^n}j(u(\lambda x),(Du)(\lambda x))+o(1) \\
&=\lambda^{p-n}\int_{\R^n}j(u,Du)+o(1),\quad \text{as $h\to\infty$}.
\end{align*}
If $L_h=(\partial_j \Pi_h^i)$ and $M_h=(\partial_j (\Pi^{-1}_h)^i)$ denote the $n\times n$
Jacobian matrices of the maps $\Pi_h,\Pi_h^{-1}$ respectively, taking into account that
$\gamma=\sup\{\|L_h\|_{n\times n}:h\in\N\}<\infty$ and $\gamma'=\sup\{|{\rm det}\,M_h|:h\in\N\}<\infty$ and using
the growth condition \eqref{growth1} on $j$, we get
\begin{align*}
I_2 &=\int_{\R^n}j(u(\Pi_h(x)),Du(\Pi_h(x)))\chi_{\{h\leq |x|\leq rh\}}  \\
& \leq c_2\int_{\R^n} |Du(\Pi_h(x))|^p\chi_{\{h\leq |x|\leq rh\}} \\
& = c_2\int_{\R^n} |L_h (Du)(\Pi_h(x))|^p\chi_{\{h\leq |x|\leq rh\}} \\
& \leq c_2\gamma^p\int_{\R^n} |(Du)(\Pi_h(x))|^p\chi_{\{h\leq |x|\leq rh\}} \\
& = c_2\gamma^p\gamma'\int_{\R^n} |Du(y)|^p\chi_{\{h\leq |\Pi_h^{-1}(y)|\leq rh\}}=o(1),\quad \text{as $h\to\infty$}.
\end{align*}
Concerning the last equality notice that, as $|\Pi_h(x)|\geq\rho
|x|\geq\rho h$, where $\rho$ is the constant appearing in
definition \eqref{defLambda}, the condition $|\Pi^{-1}_h(y)|\geq
h$ implies $|y|\geq \rho h$ and hence the integrand goes to zero
pointwise. Finally, of course, we have
\begin{align*}
I_3 &=\int_{\R^n}j(u,Du)\chi_{\{|x|>rh\}}=o(1),\quad \text{as $h\to\infty$}.
\end{align*}
In conclusion, we get
\begin{equation*}
\int_{\R^n}j(u+\phi_h,Du+D\phi_h)=\lambda^{p-n}\int_{\R^n}j(u,Du)+o(1),\quad \text{as $h\to\infty$}.
\end{equation*}
In the same way, we get
\begin{align}
\label{Fasssintt}
\int_{\R^n}F(u+\phi_h) & =\int_{\R^n}F(u(\Pi_h(x)))=\int_{\R^n}F(u(y))|{\rm det}(M_h)| \\
& =\lambda^{-n}\int_{\R^n}F(u)+o(1),\quad \text{as $h\to\infty$}.\notag
\end{align}
Finally, collecting the previous formulas, we reach the inequality
\begin{equation*}
T+(\lambda^{p-n}-1)\int_{\R^n} j(u,Du)\geq T\left(
1+(\lambda^{-n}-1)\int_{\R^n}F(u)  \right)^{\frac{n-p}{n}}\!\!\!\!
\end{equation*}
which holds for every $\lambda$ sufficiently close to $1$.
Choosing $\lambda=1+\omega$ and $\lambda=1-\omega$ with $\omega>0$ small
and then letting $\omega\to 0^+$, we conclude that
$$
\int_{\R^n} j(u,Du)=T\int_{\R^n}F(u).
$$
Since $u\not\equiv 0$, it follows that $\int_{\R^n}F(u)>0$, so
that plugging $w=u$ into \eqref{rescdis} one entails
$\int_{\R^n}F(u)\geq 1$. On the other hand, inequality
\eqref{lowersem} yields $\int_{\R^n} j(u,Du)\leq T$. This, of
course, forces
\begin{equation}
\label{conclusion}
T=\int_{\R^n} j(u,Du),\qquad \int_{\R^n}F(u)=1,
\end{equation}
which concludes the proof.
\end{proof}

\begin{remark}
\label{subcritic}
In the proof of Lemma~\ref{lemma1}, in order to show that the (minimizing) sequence $(u_h)$ is
strongly convergent in $D^{1,p}(\Omega)$ for any bounded domain $\Omega$ of $\R^n$,
we have exploited the sub-criticality assumption \eqref{inftycondit} on $f$, which is stronger than the corresponding
assumption (2.6) in \cite{brezislieb} on $F$, that is
$$
\lim_{s\to\infty}\frac{F(s)}{|s|^{p^*}}=0.
$$
In \cite{brezislieb}, due to the particular structure of $j$,
namely the model case $j(s,\xi)=\frac{1}{2}|\xi|^2$, to conclude
the proof the weak convergence of $(u_h)$ to $u$ in $D^{1,2}(\R^n)$ turns out to be
sufficient, while to cover the general case $j(s,\xi)$ the local
convergence  seems to be
necessary to handle the perturbation argument devised at the end of Lemma~\ref{lemma1}.
We point out that, also in \cite{gazfer}, the authors assume condition \eqref{inftycondit} on $f$, although they are
allowed to take a spherically symmetric minimizing sequence which provides compactness.
\end{remark}

\begin{lemma}
\label{lemma2}
Assume \eqref{zerocondit}-\eqref{posit}. Then condition $\boldsymbol{(C2)}$ holds.
\end{lemma}
\begin{proof}
Let $u\in D^{1,p}(\R^n)$ be a minimizer for problem $(P_1)$. Then
the sequence $u_h=u$ is minimizing for $(P_1)$. By Proposition
\ref{ekelVP} we can find a sequence $(v_h)\subset D^{1,p}(\R^n)$
such that $\|v_h-u\|_{D^{1,p}}=o(1)$ and $|dJ|_{{\mathcal
C}}|(v_h)\to 0$ as $h\to\infty$. Hence, by Lemma~\ref{quasieuler}
there exists a sequence $(\mu_h)\subset\R$ such that
\begin{equation}
\label{quasi-solBis} J'(v_h)(\varphi)=\mu_h
V'(v_h)(\varphi)+\langle \eta_h,\varphi\rangle,\quad\text{for all
$\varphi\in C^{\infty}_c(\R^n)$}
\end{equation}
where $\eta_h$ converges strongly to $0$ in $D^*$. As in the proof
of Lemma \ref{lemma1} we can assume that $\mu_h \to \mu$ and since
$v_h\to u$ in $D^{1,p}(\R^n)$, obviously
$$
J'(u)(\varphi)=\mu V'(u)(\varphi),\quad\text{for all $\varphi\in
C^{\infty}_c(\R^n)$}.
$$
Namely equation \eqref{autonomGC} is satisfied in the sense of
distributions for some $\mu\in\R$. By means of
assumptions~\eqref{growth1}-\eqref{growth2} and~\eqref{posit}, a
standard argument yields $u\in L^\infty_{{\rm loc}}(\R^n)$ (see,
e.g.,~\cite[Theorem 1 and Remark at p.261]{serrin}). By the
regularity results contained in \cite{benedetto,tolks}, it follows
that $u\in C^{1,\beta}_{{\rm loc}}(\R^n)$, for some $0<\beta<1$.
\end{proof}

Let $\varphi\in L^{\infty}_{{\rm loc}}(\R^n)$ and  let $L(s,\xi):\R\times\R^n\to\R$ be a function of class $C^1$ in
$s$ and $\xi$ such that, for any $s\in\R$, the map $\{\xi\mapsto L(s,\xi)\}$ is strictly convex.
We recall, in the autonomous setting, a Pucci-Serrin variational
identity for locally Lipschitz continuous solutions of a general class of equations,
recently obtained in \cite{dms}.

\begin{lemma}
\label{step1}
Let $u:\R^n\to\R$ be a locally Lipschitz solution of
\begin{equation*}
- \dvg \!\left(L_{\xi}(u,Du)\right) +
L_s(u,Du)=\varphi\quad
\text{in ${\mathcal D}'(\R^n)$}.
\end{equation*}
Then
\begin{equation}
\label{prima}
\sum_{i,j=1}^n\int_{\R^n} D_ih^j
D_{\xi_i}L(u,Du) D_ju
- \int_{\R^n}(\dvg h) \, L(u,Du)=
\int_{\R^n} (h\cdot Du)\varphi
\end{equation}
for every $h\in C^1_c(\R^n,\R^n)$.
\end{lemma}
\smallskip

\begin{remark}
The classical Pucci-Serrin identity \cite{ps} is not applicable here, since it
requires the $C^2$ regularity of the solutions while in our degenerate setting (for $p\neq 2$)
the maximal allowed regularity is $C^{1,\beta}_{{\rm loc}}$ (see \cite{benedetto,tolks}).
\end{remark}

\begin{lemma}
\label{lemma3}
Assume \eqref{zerocondit}-\eqref{posit}. Then condition $\boldsymbol{(C3)}$ holds.
\end{lemma}
\begin{proof}
Let $u\in D^{1,p}(\R^n)$ be any solution of equation \eqref{autonomGC}. In light of conditions~\eqref{growth1}, \eqref{growth2} and~\eqref{posit},
as we observed in the proof of Lemma~\ref{lemma2}, it follows that $u\in C^{1,\beta}_{{\rm loc}}(\R^n)$ for some $0<\beta<1$.
Then, since $\{\xi\mapsto j(s,\xi)\}$ is strictly convex,
we can use Lemma~\ref{step1} by choosing in \eqref{prima}
$\varphi=0$ and
\begin{align}
\label{sceltal}
L(s,\xi)&=j(s,\xi)-\mu F(s),\qquad
\text{for all $s\in\R^+$ and $\xi\in\R^n$},\\
h(x)&=h_k(x)=T\left(\frac{x}{k}\right)x,\qquad
\text{for all $x\in\R^n$ and $k\geq 1$},\notag
\end{align}
being $T\in C_c^1(\R^n)$ such that
$T(x)=1$ if $|x|\leq 1$ and $T(x)=0$ if $|x|\geq 2$.
In particular, for every $k$ we have
that $h_k\in C^1_c(\R^n,\R^n)$ and
\begin{align*}
D_ih^j_k(x)&=D_iT\left(\frac{x}{k}\right)\frac{x_j}{k}+
T\left(\frac{x}{k}\right)\delta_{ij},\qquad
\text{for all $x\in\R^n$, $i,j=1,\dots, n$,}\\
({\rm div}\,h_k)(x)&=DT\left(\frac{x}{k}\right)\cdot\frac{x}{k}+
nT\left(\frac{x}{k}\right),\qquad
\text{for all $x\in\R^n$}.
\end{align*}
Then it follows by identity \eqref{prima} that
\begin{multline*}
\sum_{i,j=1}^n\int_{\R^n}
D_iT\left(\frac{x}{k}\right)\frac{x_j}{k}D_ju
D_{\xi_i}L(u,Du)
+\int_{\R^n}T\left(\frac{x}{k}\right)D_{\xi}L(u,Du)\cdot Du+\\
-\int_{\R^n}DT\left(\frac{x}{k}\right)\cdot\frac{x}{k}\, L(u,Du)-
\int_{\R^n}nT\left(\frac{x}{k}\right)L(u,Du)=0,
\end{multline*}
for every $k\geq 1$. Since there exists $C>0$ with
\begin{equation*}
\big|D_iT\left(\frac{x}{k}\right)\frac{x_j}{k}\big|\leq C\quad
\text{for every $x\in\R^n$, $k\geq 1$ and $i,j=1,\dots,n$,}
\end{equation*}
by the Dominated Convergence Theorem (recall that by \eqref{growth1} and the $p$-homogeneity
of $\{\xi\mapsto j(s,\xi)\}$, of course one has $L(u,Du),D_{\xi}L(u,Du)\cdot Du\in L^1(\R^n)$),
letting $k\to\infty$, we conclude that
\begin{gather*}
\int_{\R^n}\Big[n L(u,Du)-
D_{\xi}L(u,Du)\cdot Du\Big]=0,
\end{gather*}
namely, by \eqref{sceltal} and, again, the
$p$-homogeneity of $\{\xi\mapsto j(s,\xi)\}$,
\begin{equation}
\label{pucciser}
(n-p)\int_{\R^n}j(u,Du)=\mu n\int_{\R^n}F(u),
\end{equation}
namely $(n-p)J(u)=\mu nV(u)$, proving that condition $\boldsymbol{(C3)}$ is fulfilled.
\end{proof}
\smallskip

\noindent{\bf Proof of Theorem \ref{MainThm}.} From Lemmas
\ref{lemma1}, \ref{lemma2} and \ref{lemma3} we see that the
conditions {\bf (C1)-(C3)} hold. The conclusion follows directly
from Proposition \ref{abstract}. \hfill $\Box$

\begin{remark}
\label{observ-positivity} In light of formula \eqref{pucciser} and
the positivity of $j$, it holds $\int_{\R^n}F(u)>0$ as soon as $u$
is a nontrivial solution of \eqref{autonomGC}.
\end{remark}

\begin{remark}
\label{rem27} Assumption~\eqref{posit} was already considered
e.g.\ in \cite{bbm,bmp,canino,cd,giacsqua,squastop,squastoul}. We
exploited it in order to get existence (in $\boldsymbol{(C_1)}$),
regularity (in $\boldsymbol{(C_2)}$) and hence also for the
Pucci-Serrin identity (in $\boldsymbol{(C_3)}$), and it seems hard
to drop, mainly concerning the boundedness (and hence $C^1$ regularity) issue of
solutions. In fact, in lack of~\eqref{posit} some problems may
occur, already in the case of bounded domains and $p=2$. For
instance, as shown by J.\ Frehse in \cite{fresh}, if $B(0,1)$ is
the unit ball in $\R^n$ centered at zero with $n\geq 3$,
$$
j(x,s,\xi)=\Big(1 + \frac{1}{|x|^{12(n-2)}\,e^s+1}\Big)|\xi|^2
$$
and $f(s)=0$, then $u(x)=-12(n-2)\log|x|$ is a weak solution to
the corresponding Euler equation with $u=0$ on $\partial
B(0,1)$. In particular $u\not\in L^{\infty}(B(0,1))$ although $j$ is very
regular. It is immediate to check that $j_s(x,s,\xi)s\leq 0$ for
any $s\geq 0$, so ~\eqref{posit} fails.
Although this counterexample involves an $x$-dependent Lagrangian (while we deal
with autonomous problems) these pathologies in regularity are related to
the $s$-dependence in the Lagrangian $j$.
\end{remark}
\bigskip


\section{The case $p=n$}

 We consider
now the following constrained minimization problem
\begin{equation}
\label{constminim2} \tag{$P_0$} \text{minimize $J(u)$ for $u\neq
0$ subject to the constraint $V(u)=0$}.
\end{equation}
More precisely, let us set
\begin{equation}
\label{defX2} X_0=\big\{u\in D^{1,n}(\R^n):\, F(u)\in
L^1(\R^n)\big\},
\end{equation}
and
$$
T_0=\inf_{{\mathcal C}_0}J,\qquad {\mathcal C}_0=\big\{u\in
X_0:\,\,u\neq 0,\,\,V(u)=0\big\}.
$$
Consider the following conditions: \vskip6pt \noindent
$\boldsymbol{(D1)}$ $T>0$ and problem $(P_0)$ has a minimizer
$u\in X_0$; \vskip8pt \noindent $\boldsymbol{(D2)}$ any minimizer
$u\in X_0$ of $(P_0)$ is a $C^1$ solution and satisfies the
equation
\begin{equation}
\label{autonomGC2} -\dvg(j_\xi(u,Du))+j_s(u,Du)=\mu f(u)\quad
\text{in ${\mathcal D}'(\R^n)$},
\end{equation}
for some $\mu\in\R$. \vskip8pt \noindent $\boldsymbol{(D3)}$ any
solution $u\in X_0$ of equation \eqref{autonomGC2} with $\mu>0$
satisfies $V(u)=0$. \vskip6pt \noindent

From \cite{ByJeMa} we have the following

\begin{proposition} \label{abstract2}
Assume that $p=n$ and that {\bf (D1)-(D3)} hold. Then
(\ref{autonom}) admits a least energy solution and each least
energy solution has a constant sign. Moreover if  $u \in
D^{1,n}(\R^n)$ is a least energy solution such that $u(x) \to 0$
as $|x| \to \infty$ it is radially symmetric, up to a translation in
$\R^n$.
\end{proposition}

Proposition \ref{abstract2} follows directly from Propositions 4
and 6 in \cite{ByJeMa}. See also Remark \ref{adaptation}.
\medskip

Let us now show that {\bf (D1)-(D3)} hold. First we recall a
regularity result (see \cite{bb}).
\begin{lemma}
\label{bb} Let $u,v\in D^{1,n}(\R^n),$ $\eta\in L^1(\R^n)$ and
$w\in D^*(\R^n)$ with
$$
j_s(u,\nabla u)v\geq\eta,
$$
and for all $\varphi\in C^\infty_c(\R^n)$
\begin{equation*}
\langle w,\varphi\rangle=\int_{\R^n} j_\xi(u,\nabla u) \cdot D
\varphi+\int_{\R^n} j_s(u,Du)\varphi.
\end{equation*}
Then $j_s(u,\nabla u)v\in L^1(\R^n)$ and
\begin{equation*}
\langle w,v\rangle=\int_{\R^n} j_\xi(u,\nabla u)\cdot
Dv+\int_{\R^n} j_s(u,\nabla u)v.
\end{equation*}
\end{lemma}

Now we have

\begin{proposition}
\label{lemma1-2} Assume \eqref{stricconv}-\eqref{manifold2}. Then
conditions $\boldsymbol{(D1)}$-$\boldsymbol{(D3)}$ hold.
\end{proposition}
\begin{proof}
In view of \eqref{existposit2} the constraint ${\mathcal C}_0$ is
not empty (see again Step 1 at page 324 in \cite{berlions}). Let
then $(u_h)\subset {\mathcal C}_0$ be a minimizing sequence for
$J|_{{\mathcal C}_0}$. Therefore, we have
$$
\lim_{h}\int_{\R^n}j(u_h,Du_h)=T_0,\quad u_h\neq 0,\quad F(u_h)\in
L^1(\R^n),\quad \int_{\R^n}F(u_h)=0,
$$
for all $h\in\N$. Since $u_h\neq 0$ and by \eqref{zerocondit2}, it
holds
$$
\int_{\{|u_h|>\delta\}} F(u_h)=\int_{\{0\leq |u_h|\leq \delta\}}
|F(u_h)|>0,
$$
of course ${\mathcal L}^n(\{|u_h|>\delta\})>0$ for every $h\in\N$.
Then, since the map $\{u\mapsto J(u)\}$ is invariant under scaling
on $D^{1,n}(\R^n)$, it is readily seen that there exists
$\varrho>0$ such that
\begin{equation} \label{nonvanishing}
{\mathcal L}^n(\{|u_h|>\delta\})\geq\varrho,\qquad\text{for all
$h\in\N$}.
\end{equation}
Arguing as in \cite[Lemma 3.1]{brezislieb}, we have that
\begin{equation}
\label{summabcond} \sup_{h\in\N}\int_{\{|u_h|>\delta\}}
|u_h|^{r}<\infty,\qquad\text{for all $r>1$}.
\end{equation}
Then the sequence $(u_h)$ is bounded in $L^q(\Omega)$ for any
bounded domain $\Omega\subset\R^n$ and, after extracting a
subsequence still denoted by $(u_h)$, we have $u_h\to u$ in
$L^{q}(\Omega)$ for all $q \geq 1$, $Du_h\rightharpoonup Du$ in
$L^{n}(\R^n)$, and $u_h(x)\to u(x)$ a.e. $x\in\R^n$. As $j(s,\xi)$
is positive, convex in the second argument, $u_h\to u$ in
$L^1_{{\rm loc}}(\R^n)$ and $Du_h\rightharpoonup Du$ in $L^1_{{\rm
loc}}(\R^n)$ by lower semi-continuity it follows
\begin{equation}
\label{lowersem2} \int_{\R^n}j(u,Du)\leq
\liminf_{h}\int_{\R^n}j(u_h,Du_h)=T_0.
\end{equation}

Now from (\ref{nonvanishing}), as in the proof of Lemma
\ref{lemma1} we get that, after a shift, the weak limit of $(u_h)$
is nontrivial, that is $u\neq 0$. Notice also that, in view of
\eqref{zerocondit2}, \eqref{inftycondit2} and the bound furnished
by \eqref{summabcond} we get, for any $h\in\N$,
$$
\int_{\R^n} F_-(u_h)=\int_{\R^n} F_+(u_h)=\int_{\{|u_h|>\delta\}}
F_+(u_h)\leq C\int_{\{|u_h|>\delta\}} |u_h|^{q}\leq C,
$$
where $C$ is a generic positive constant. In particular, by Fatou's
lemma, it follows $F\in L^1(\R^n)$. We have proved thus $u\in
X_0$. Arguing as in Lemma~\ref{lemma1}, up to substituting
$(u_h)\subset {\mathcal C}_0$ with a new minimizing sequence
$(v_h)\subset {\mathcal C}_0$, we may assume that $|dJ|_{{\mathcal
C}_0}|(u_h)\leq\eps_h$, with $\eps_h\to 0$ as $h\to\infty$. By
Lemma~\ref{quasieuler} there exists a sequence $(\mu_h)\subset\R$
such that
\begin{equation}
\label{quasi-sol2} J'(u_h)(v)=\mu_h V'(u_h)(v)+\langle
\eta_h,v\rangle,\qquad\text{for all $h\in\N$ and $v\in
C^{\infty}_c(\R^n)$,}
\end{equation}
where $\eta_h$ strongly converges to $0$ in $D^*$ as $h\to\infty$.
As in Lemma~\ref{lemma1}, it can be proved that $(\mu_h)$ is
bounded (and hence it converges to some value $\mu\in\R$). Now
since \eqref{quasi-sol2} hold and $u_h \rightharpoonup u$ in
$D^{1,n}(\R^n)$ using  the classical convergence result of Murat
(see Theorem 2.1 of \cite{bm}) we get that $Du_h(x) \to Du(x)$
a.e. $x \in \R^n$. At this point it follows easily that
\eqref{autonomGC2} is satisfied (see e.g.\ \cite[Theorem 3.4]{squastoul} for details).

Let us now prove that, actually, $\mu\neq 0$. If, by contradiction, it
was $\mu=0$, then we would have
$$
\int_{\R^n}j_\xi(u,Du)\cdot
Dv+\int_{\R^n}j_s(u,Du)v=0,\qquad\text{for all $v\in
C^{\infty}_c(\R^n)$}.
$$
Let now $\zeta:\R\to\R$ be the map defined by
\begin{equation}
\label{zeta} \zeta(s)=
\begin{cases}
M|s| & \text{if $|s|\leq R$} \\
MR & \text{if $|s|\geq R$},
\end{cases}
\end{equation}
being $R>0$ the constant defined in \eqref{posit} and $M$ a
positive number (which exists by combining the growths conditions
\eqref{growth1}-\eqref{growth2}) such that
\begin{gather}
\label{controlBB} \left|j_s(s,\xi)\right|\leq nMj(s,\xi)
\end{gather}
for $s\in\R$ and $\xi\in\R^n$. Notice that, by combining
\eqref{posit} and \eqref{controlBB}, we obtain
\begin{equation}
\label{poss} \big[j_s(s,\xi)+n\zeta'(s)j(s,\xi)\big]s\geq 0,\qquad
\text{for all $s\in\R$ and $\xi\in\R^n$}.
\end{equation}
Taking into account \eqref{posit}, by Lemma \ref{bb} we are
allowed to choose $v=ue^{\zeta(u)}$ and hence
$$
\int_{\R^n}ne^{\zeta(u)}
j(u,Du)+\int_{\R^n}e^{\zeta(u)}\big[j_s(u,Du)+n\zeta'(u)j(u,Du)\big]u=0.
$$
 Then, by \eqref{poss} and \eqref{growth1} we get
$$
nc_1\int_{\R^n}|Du|^n\leq 0,
$$
so that $u=0$, which is not possible. Hence $\mu\neq 0$.  Arguing
as in Lemma \ref{lemma3} the Pucci-Serrin identity follows,
namely, as $p=n$
\begin{equation*}
\label{pucciser2}
\int_{\R^n}F(u)=\frac{n-p}{\mu
n}\int_{\R^n}j(u,Du)=0.
\end{equation*}
The same conclusion obviously hold for any solution $u \in X_0$
for (\ref{autonomGC2}) with $\mu >0$, this shows that {\bf (D3)}
hold. Now since $u\in X_0$ and $\int_{\R^n}F(u)=0$, we have
$u\in{\mathcal C}_0$, so that by \eqref{lowersem2}
\begin{equation*}
\int_{\R^n}j(u,Du)=T_0.
\end{equation*}
As in Lemma~\ref{lemma2}, one can prove that any minimizer is
$C^1$ and satisfies the Euler-Lagrange equation
\eqref{autonomGC2}, which concludes the proof.
\end{proof}
\smallskip

\begin{remark}
The check that {\bf (D1)} holds is actually simpler than in the case of {\bf (C1)}.
In particular, to check {\bf (D1)} we do not need to use any kind of strong local
convergence, as for the case $1<p<n$, using classical convergence results due to Murat
suffices. Observe also that, in the case $p=n$, if we have a nontrivial function $v\in X_0$ which is
a solution to the problem
\begin{equation}
\label{prob-v} -\dvg(j_\xi(v,Dv))+j_s(v,Dv)= f(v)\quad \text{in ${\mathcal D}'(\R^n)$},
\end{equation}
then, by the Pucci-Serrin identity it follows that
$\int_{\R^n}F(v)=0$, so that $v\in{\mathcal C}_0$ and hence
$$
I(v)=\int_{\R^n}j(v,Dv)\geq T_0=\int_{\R^n}j(u,Du)=I(u),
$$
proving that $u$ is, automatically, a least energy solution of \eqref{prob-v}.
\end{remark}

\bigskip

\end{document}